\documentclass[12pt,a4paper]{article}
\usepackage{amsmath,amssymb,amsthm}

\def\R{\mathbb{R}}

\def\K{\mathcal{K}}

\def\fd{\mathfrak{d}}

\def\ssm{\smallsetminus}

\def\restrictedto%
{\mathchoice%
    {\!\upharpoonright\!}
    {\!\upharpoonright\!}
    {\upharpoonright}
    {\upharpoonright}
}

\def\st{\mathchoice{:}{:}{\,:\,}{\,:\,}}

\def\size#1{\lvert#1\rvert}

\def\domn{\leq^*}

\newcommand{\splitnodes}[2][{}]%
	{\operatorname{\textup{\textsf{split}}}_{#1}({#2})}

\def\abs#1{\left\lvert{#1}\right\rvert}
\def\cR{\mathcal{R}}
\def\tfrac#1#2{\textstyle\frac{#1}{#2}}
\def\cl{\operatorname{cl}}	
\def\interior{\operatorname{int}}	
\def\Ball{\operatorname{B}}
\def\diam{\operatorname{diam}}
\newcommand{\deriv}[2][{1}]{{#2}^{({#1})}}	
\def\higson#1#2{\overline{#1}^{#2}}
\def\smirnov#1#2{u_{#2}{#1}}
\def\stonecech#1{\beta{#1}}
\def\PM{\operatorname{PM}}
\def\Metric{\operatorname{M}}
\def\cptequiv{\simeq}

\def\fsa{\mathfrak{sa}}

\def\fha{\mathfrak{ha}}

\def\cptsep#1#2#3{{#1}\parallel{#2}\,\,\,(#3)}
\def\cptnonsep#1#2#3{{#1}\not\,\parallel{#2}\,\,\,(#3)}

\theoremstyle{plain}
\newtheorem{thm}{Theorem}[section]
\newtheorem{lem}[thm]{Lemma}
\newtheorem{cor}[thm]{Corollary}
\newtheorem{prop}[thm]{Proposition}

\newtheorem{claim}{Claim}

\theoremstyle{definition}
\newtheorem{defn}[thm]{Definition}
\theoremstyle{remark}

\begin{document}

\title{How many miles to $\beta X$?
	--- $\mathfrak{d}$ miles, or just one foot}
\author{Masaru Kada%
    \thanks{%
	Supported by 
	Grant-in-Aid for Young Scientists (B) 14740058, MEXT.} 
	\and
	Kazuo Tomoyasu%
	\thanks{%
	Supported by 
	Grant-in-Aid for Young Scientists (B) 14740057, MEXT.} 
	\and 
	Yasuo Yoshinobu%
	\thanks{%
	Supported by
	Grant-in-Aid for Scientific Research (C)(2) 15540115, JSPS.}
	}
\date{}
\maketitle

\begin{abstract}
It is known that 
the Stone--\v{C}ech compactification $\stonecech{X}$ 
of a 
metrizable space $X$ 
is approximated by the collection of Smirnov compactifications of $X$ 
for all compatible metrics on $X$. 
If we confine ourselves to locally compact separable metrizable spaces, 
the corresponding statement holds for Higson compactifications. 
We investigate the smallest cardinality 
of a set $D$ of compatible metrics 
on $X$ such that 
$\stonecech{X}$ is approximated by Smirnov or Higson compactifications 
for all metrics in $D$. 
We prove that 
it is either the dominating number or 1 
for a locally compact separable metrizable space. 

\par
\vspace{12pt}
 \par
\vspace{12pt}
\end{abstract}

\section{Introduction}

A \emph{compactification} 
of a completely regular Hausdorff 
space $X$ 
is a compact Hausdorff space 
which contains $X$ as a dense subspace. 
For compactifications $\alpha X$ and $\gamma X$ 
of $X$, 
we write $\alpha X\leq \gamma X$ 
if there is a continuous surjection $f\st \gamma X\to\alpha X$ 
such that 
$f\restrictedto X$ is the identity map on $X$.
If such an $f$ can be chosen to be a homeomorphism, 
we write $\alpha X\cptequiv\gamma X$. 
Let $\K(X)$ denote the class of compactifications of $X$. 
When 
we identify $\cptequiv$-equivalent compactifications, 
we may regard $\K(X)$ as a set, 
and the order structure $(\K(X),{\leq})$ 
is a complete upper semilattice 
whose largest element is the Stone--\v{C}ech compactification $\stonecech{X}$.

Let $C^*(X)$ denote 
the ring of bounded continuous functions 
from $X$ to $\R$. 
A subring $R$ of $C^*(X)$ is called \emph{regular} 
if $R$ is closed in the sense of uniform norm topology,
contains all constant functions, 
and generates the topology of $X$.
Let $\cR(X)$ denote the class of regular subrings of $C^*(X)$.
Then it is known that $(\K(X), {\leq})$ is isomorphic to
$(\cR(X), {\subseteq})$, by mapping each $\alpha X\in\K(X)$
to the set of bounded continuous functions from $X$ to $\R$ 
which are continuously extended over $\alpha X$
(cf.\ \cite[Theorem 3.7]{BaY:ring}, \cite[Theorem 2.5]{Chan:cpt}).
In particular, 
$\beta X$
corresponds to the whole $C^*(X)$.
(See \cite{Chan:cpt,GJ:rings} for more details.)

For a compactification $\alpha X$ of $X$ 
and two closed subsets $A,B$ of $X$, 
we write $\cptsep{A}{B}{\alpha X}$ 
if $\cl_{\alpha X}A\cap\cl_{\alpha X}B=\emptyset$, 
and otherwise 
$\cptnonsep{A}{B}{\alpha X}$. 

For a metric space $(X,d)$, 
$U^*_{d}(X)$ denotes the set of all bounded uniformly continuous functions 
from $(X,d)$ to $\R$.
$U^*_{d}(X)$ is a regular subring of $C^*(X)$. 
The \emph{Smirnov compactification $\smirnov{X}{d}$ of\/ $(X,d)$} 
is the unique compactification associated with the subring $U^*_{d}(X)$.
For disjoint closed subsets $A, B$ of $X$,
$\cptsep{A}{B}{\smirnov{X}{d}}$ if and only if $d(A,B)>0$
\cite[Theorem 2.5]{Wo:unifcpt}.

The following theorem tells us that 
we can approximate 
the Stone--\v{C}ech compactification of a metrizable space 
by the collection of all Smirnov compactifications. 
Let $\Metric(X)$ denote the set of all metrics on $X$ 
which are compatible with the topology on $X$. 

\begin{thm}\label{thm:smirnovapprox}
\textup{\cite[Theorem 2.11]{Wo:unifcpt}}
For a noncompact metrizable space $X$, 
we have 
$\stonecech{X}\cptequiv\sup\{\smirnov{X}{d}\st d\in\Metric(X)\}$ 
\textup{(}%
the supremum
is taken in the lattice $(\mathcal{K}(X),{\leq})$\textup{)}. 
\end{thm}

Now we define 
the following cardinal function. 

\begin{defn}
\textup{\cite[Definition~2.2]{KTY:babylon}}
For a noncompact 
metrizable space $X$, 
let 
$\fsa(X)=\min\{\size{D}\st
	D\subseteq\Metric(X)\text{ and }
	\stonecech{X}\cptequiv\sup\{\smirnov{X}{d}\st d\in D\}\}$. 
\end{defn}

For a metrizable space $X$, 
a metric $d$ on $X$ is called \emph{proper} 
if each $d$-bounded set has compact closure. 
A \emph{proper metric space} means 
a metric space whose metric is proper. 

For a function $f$  and a subset $A$ of the domain of $f$, 
$f''A$ denotes the image of $A$ by $f$. 

Let $(X,d)$ be a proper metric space and $(Y,\rho)$ a metric space. 
We say 
a function 
$f$ from $X$ to $Y$ is \emph{slowly oscillating} 
if it satisfies the following condition: 
\begin{equation*}
\forall r>0\,
\forall\varepsilon>0\, 
\exists K\text{ a compact subset of }X\, 	
\forall x\in X\ssm K\,(\diam_\rho(f''\Ball_d(x,r))<\varepsilon). 
\end{equation*}
For a proper metric space $(X,d)$, 
let $C^*_d(X)$ be 
the set of all bounded continuous slowly oscillating functions 
from $(X,d)$ to $\R$.
$C^*_d(X)$ is a regular subring of $C^*(X)$.
The \emph{Higson compactification $\higson{X}{d}$ of $(X,d)$} 
is the unique compactification associated with the subring $C^*_d(X)$.
For disjoint closed subsets $A,B$ of $X$,
$\cptsep{A}{B}{\higson{X}{d}}$
if and only if
for any $R>0$ there is a compact subset $K_R$ of $X$ 
such that 
$d(x,A)+d(x,B)>R$ holds for all $x\in X\ssm K_R$ 
\cite[Proposition 2.3]{DKU:higsoncorona}.

The following 
corresponds to Theorem \ref{thm:smirnovapprox} 
for Higson compactifications. 
Note that 
a proper metric space 
is locally compact and separable. 
Let $\PM(X)$ be the set of all proper metrics 
compatible with the topology of $X$. 

\begin{thm}\label{thm:higsonapprox}
\textup{\cite[Theorem 3.2]{KT:approx}}
For a noncompact locally compact separable metrizable space $X$, 
we have 
$\stonecech{X}\cptequiv\sup\{\higson{X}{d}\st d\in\PM(X)\}$. 
\end{thm}

So we consider the following cardinal function. 

\begin{defn}
\textup{\cite[Definition~6.2]{KTY:babylon}}
For a noncompact locally compact separable metrizable space $X$, 
let 
$\fha(X)=\min\{\size{D}\st
	D\subseteq\PM(X)\text{ and }
	\stonecech{X}\cptequiv\sup\{\higson{X}{d}\st d\in D\}\}$. 
\end{defn}

We have $\fsa(X)\leq\fha(X)$ 
for each locally compact separable metrizable space $X$ 
\cite[Lemma~6.3]{KTY:babylon}. 

For $f,g\in\omega^\omega$, 
we say $f\domn g$ if for all but finitely many $n<\omega$ 
we have $f(n)\leq g(n)$. 
The \emph{dominating number} $\fd$ 
is the smallest size of a subset of $\omega^\omega$ 
which is cofinal in $\omega^\omega$ with respect to $\domn$. 

In Section \ref{sec:dor1} 
we will show that, 
for a locally compact separable metrizable space $X$, 
either $\fsa(X)=\fha(X)=\fd$ or $\fsa(X)=\fha(X)=1$ holds.   
In Section \ref{sec:morethand} 
we will give 
an example of a nonseparable metrizable space $X$ 
for which $\fsa(X)>\fd$ holds.

\section{Dichotomy for locally compact separable spaces}\label{sec:dor1}

It is easily seen that $\fsa(\omega)=\fha(\omega)=1$. 
In fact, 
the following two theorems give 
equivalent conditions respectively for $\fsa(X)=1$ and $\fha(X)=1$. 

For a space $X$, 
$\deriv{X}$ denotes 
the first Cantor--Bendixson derivative of $X$, 
that is, 
the subspace of $X$ which consists of 
nonisolated points of $X$. 

\begin{thm}\label{thm:saone}
\textup{\cite[Corollary~3.5]{Wo:unifcpt}}
For a metrizable space $X$, the following conditions are equivalent.
\begin{enumerate}
\item There is a metric $d\in\Metric(X)$ 
	for which $\smirnov{X}{d}\cptequiv\stonecech{X}$ holds. 
\item 
	$\deriv{X}$ is compact. 
\end{enumerate}
\end{thm}

\begin{thm}\label{thm:haone}
\textup{\cite[Proposition~2.6]{KT:approx}}
For a locally compact separable metrizable space $X$, 
the following conditions are equivalent.
\begin{enumerate}
\item There is a proper metric $d\in\PM(X)$ 
	for which $\higson{X}{d}\cptequiv\stonecech{X}$ holds. 
\item 
	$\deriv{X}$ is compact. 
\end{enumerate}
\end{thm}

In the paper \cite{KTY:babylon} 
we proved the following proposition. 

\begin{prop}\label{prop:halfopen}
\textup{\cite[Examples~2.3 and 6.4]{KTY:babylon}}
$\fsa([0,\infty))=\fha([0,\infty))=\fd$. 
\end{prop}

In this section we prove that, 
assuming that $X$ is locally compact and separable, 
$\fsa(X)=\fha(X)=\fd$ 
unless 
$\deriv{X}$ 
is compact. 
In particular, 
since $\fha(X)$ is defined only 
when $X$ is locally compact and separable, 
$\fha(X)$ is either $\fd$ or $1$ 
when it is defined. 

We will use the following two lemmas. 

\begin{lem}\label{lem:charcptsep}
\textup{\cite[Lemma~1.1]{KTY:babylon}}
For a compactification $\alpha X$ of a space $X$ 
and closed subsets $A,B$ of $X$, 
the following 
conditions 
are equivalent\/\textup{:}
\begin{enumerate}
\item $\cptsep{A}{B}{\alpha X}$. 
\item There are $g\in C^*(X)$ 
	and $a,b\in\R$ 
	such that 
	$a>b$, 
	$g(x)\geq a$ for all $x\in A$, 
	$g(x)\leq b$ for all $x\in B$ 
	and $g$ is continuously extended over $\alpha X$. 
\end{enumerate}
\end{lem}

Note that, 
for a normal space $X$, 
$\alpha X\cptequiv\stonecech{X}$ 
if and only if 
$\cptsep{A}{B}{\alpha X}$ for any disjoint closed subsets $A,B$ of $X$. 

\begin{lem}\label{lem:cptsepcompact}
\textup{\cite[Lemma~1.2]{KTY:babylon}}
Suppose that $\mathcal{C}$ is a set of compactifications of a space $X$. 
For closed sets $A,B$ of $X$, 
the following 
conditions 
are equivalent: 
\begin{enumerate}
\item $\cptsep{A}{B}{\sup\mathcal{C}}$. 
\item $\cptsep{A}{B}{\sup\mathcal{F}}$ 
	for some nonempty finite subset $\mathcal{F}$ of\/ $\mathcal{C}$. 
\end{enumerate}
\end{lem}

Since $\fsa(X)\leq\fha(X)$ holds 
if both are defined, 
it suffices to show that $\fsa(X)\geq\fd$ and $\fha(X)\leq\fd$. 

First we show that $\fsa(X)\geq\fd$ unless $\fsa(X)=1$. 
This holds for all metrizable spaces. 

\begin{lem}\label{lem:lower}
Let $X$ be a metrizable space. 
If 
$\deriv{X}$ 
is not compact, 
then $\fsa(X)\geq\fd$. 
\end{lem}

\begin{proof}
Since 
$\deriv{X}$ 
is 
not 
compact, 
there is a countable subset $A$ of 
$\deriv{X}$ 
which has no accumulating point in $X$. 
Note that $A$ is closed in $X$. 
Enumerate $A$ as $\{a_n\st n<\omega\}$. 

\begin{claim}
There are a neighborhood $U_n$ of $a_n$ 
and a sequence 
$\langle b_{n,i}\st i<\omega\rangle$ in $U_n\ssm\{a_n\}$ 
for $n<\omega$ 
such that, 
\begin{enumerate}
\item for each $n<\omega$, 
	$\langle b_{n,i}\st i<\omega\rangle$ converges to $a_n$, 
\item if $n<m<\omega$ then $U_n\cap U_m=\emptyset$, and 
\item for any $f\in\omega^\omega$, 
	the set $B_f=\{b_{n,f(n)}\st n<\omega\}$ 
	has no accumulating point. 
\end{enumerate}
\end{claim}

\begin{proof}
Fix a metric $\rho\in\Metric(X)$. 
For each $n<\omega$, 
let 
$\delta_n=\frac{1}{3}\cdot\rho(a_n,A\ssm\{a_n\})$. 
By the choice of $A$, we have $\delta_n>0$. 
Let $U_n=\Ball_\rho(a_n,\delta_n)$. 
Then $n\neq m$ implies $U_n\cap U_m=\emptyset$. 
Since $a_n$ is not isolated in $X$, 
we can choose a sequence $\langle b_{n,i}\st i<\omega\rangle$ 
in $U_n\ssm\{a_n\}$ which converges to $a_n$. 
Fix an arbitrary $f\in\omega^\omega$. 
By the choice of $\delta_n$'s, 
if $B_f$ accumulates to a point, 
then $A$ must accumulate to the same point. 
Hence $B_f$ has no accumulating point. 
\end{proof}

Fix $\kappa<\fd$ and a set $D\subseteq\Metric(X)$ of size $\kappa$. 
We show that $\stonecech{X}\not\cptequiv\sup\{\smirnov{X}{d}\st d\in D\}$. 

For each $d\in D$, 
define 
a function 
$g_d\in\omega^\omega$ by letting 
\[
	g_d(n)=\min\left\{m<\omega\st
		\forall i\geq m\,
		\left(d(a_n,b_{n,i})<\tfrac{1}{n+1}\right)\right\}
\]
for $n<\omega$. 
For each nonempty finite subset $F$ of $D$, 
let $g_F=\max\{g_f\st f\in F\}$ 
(where $\max$ is the pointwise maximum). 
Since 
$\size{[D]^{<\omega}}=\size{D}=\kappa<\fd$, 
there is an $f\in\omega^\omega$ 
which satisfies 
$f\not\domn g_F$ 
for every nonempty finite subset $F$ of $D$. 

Let $B=B_f=\{b_{n,f(n)}\st n<\omega\}$. 
Then $B$ is closed and disjoint from $A$. 

For an arbitrary nonempty finite subset $F$ of $D$, 
the set $I_F=\{n<\omega\st g_F(n)<f(n)\}$ 
is infinite. 
Let $C=\cl\langle\bigcup\{U^*_d(X)\st d\in F\}\rangle$. 
Then $C$ is the closed subring of $C^*(X)$ 
associated with $\sup\{\smirnov{X}{d}\st d\in F\}$. 
By the definition of $g_F$, 
each $n\in I_F$ satisfies 
$d(a_n,b_{n,f(n)})<\frac{1}{n+1}$ for all $d\in F$. 
If $\psi\in\bigcup\{U^*_d(X)\st d\in F\}$, 
then the sequence 
$\langle\psi(a_n)-\psi(b_{n,f(n)})\st n\in I_F\rangle$ 
converges to $0$. 
So for all $\varphi\in C$, 
$\langle\varphi(a_n)-\varphi(b_{n,f(n)})\st n\in I_F\rangle$ 
converges to $0$. 
This means that 
there are no $\varphi\in C$ and $a,b\in\R$ 
such that 
$a>b$, 
$\varphi(x)\geq a$ for all $x\in A$, 
and 
$\varphi(x)\leq b$ for all $x\in B$. 
By Lemma~\ref{lem:charcptsep}, 
this means 
$\cptnonsep{A}{B}{\sup\{\smirnov{X}{d}\st d\in F\}}$. 
Since $F$ is an arbitrary nonempty finite subset of $D$ 
and by 
Lemma~\ref{lem:cptsepcompact}, 
we have 
$\cptnonsep{A}{B}{\sup\{\smirnov{X}{d}\st d\in D\}}$, 
and hence 
$\stonecech{X}\not\cptequiv\sup\{\smirnov{X}{d}\st d\in D\}$. 
\end{proof}

We turn to the proof of 
the inequality 
$\fha(X)\leq\fd$. 

For notational convenience, 
in the following lemmas and proofs, 
we let $C_n=K_n=\emptyset$ for $n=-1,-2,\ldots$

\begin{lem}\label{lem:stepfunction}
Suppose that $X$ is a normal space, 
and a sequence 
$\langle C_n\st n<\omega\rangle$ 
of closed subsets of $X$ 
satisfies $C_n\subseteq\interior C_{n+1}$ for all $n<\omega$ 
and $X=\bigcup\{C_n\st n<\omega\}$. 
Then, 
for an increasing sequence 
$\langle r_n\st n<\omega\rangle$ 
of nonnegative real numbers, 
there is a continuous function $\varphi$ from $X$ to $[0,\infty)$ 
such that, 
for each $n<\omega$ 
we have $\varphi''(C_n\ssm\interior C_{n-1})\subseteq[r_n,r_{n+1}]$. 
\end{lem}

\begin{proof}
For each $n<\omega$, 
choose a continuous function $\varphi_n$ from $X$ to $[0,r_n]$ 
so that 
${\varphi_n}''C_{n-2}=\{0\}$ 
and 
${\varphi_n}''(X\ssm\interior C_{n-1})=\{r_n\}$. 
Note that, 
if $x\in C_m$, 
then for all $n\geq m+2$ we have $\varphi_n(x)=0$. 
So we can define a continuous function $\varphi$ from $X$ to $[0,\infty)$ 
as the pointwise maximum of $\{\varphi_n\st n<\omega\}$, 
and then $\varphi$ satisfies the requirement. 
\end{proof}

Suppose that $X$ is a locally compact separable metrizable space. 
Since $X$ is $\sigma$-compact, 
there is a sequence $\langle K_n\st n<\omega\rangle$ 
of compact subsets of $X$ such that, 
for each $n<\omega$ we have $K_n\subseteq\interior K_{n+1}$, 
and $X=\bigcup\{K_n\st n<\omega\}$. 

\begin{lem}\label{lem:expand}
Let $(X,d)$ be a locally compact separable metric space, 
and 
$\langle K_n\st n<\omega\rangle$ 
a sequence of compact subsets of $X$ such that, 
for each $n<\omega$ we have $K_n\subseteq\interior K_{n+1}$, 
and $X=\bigcup\{K_n\st n<\omega\}$. 
Then, 
for each $g\in\omega^\omega$, 
there is a proper metric $d_g$ which satisfies the following:
\begin{enumerate}
\item $d_g$ is compatible with the topology of $X$. 
\item For $n<\omega$ and $x,y\in X\ssm K_{n-1}$ 
	we have $d_g(x,y)\geq g(n)\cdot d(x,y)$. 
\item For $n<\omega$ 
	we have $d_g(K_{n-1},X\ssm K_n)\geq n$. 
\end{enumerate}
\end{lem}

\setcounter{claim}{0}

\begin{proof}
Let $R_n=\max\{n,\diam_d(K_n)\}$ for each $n<\omega$, 
and let $c$ be the continuous function from $X$ to $[0,\infty)$ 
which is obtained by applying 
Lemma~\ref{lem:stepfunction} 
to 
$\langle K_n\st n<\omega\rangle$ 
and 
$\langle R_n\st n<\omega\rangle$. 

We may assume that $g$ is increasing and $g(0)\geq 1$. 
Choose an increasing continuous function $f$ 
from $[0,\infty)$ to $[1,\infty)$ 
such that 
$f(\frac{n}{2})\geq g(n)$ for all $n<\omega$. 
For $s\in [0,\infty)$, let 
\[
F(s)=\int_0^s f(t)dt.
\]
Define functions $\rho$, $\rho'_g$ from $X\times X$ to $[0,\infty)$ 
by the following:
\[
\rho(x,y)=\max\{\abs{c(x)-c(y)},d(x,y)\},
\]
\[
\rho'_g(x,y)
=	f(\max\{c(x),c(y)\})\cdot\rho(x,y).
\]
It is 
easy
to see that 
$\rho$ is a proper metric on $X$ 
and compatible with the topology on $X$. 
However, 
$\rho'_g$ is not necessarily a metric on $X$, 
because $\rho'_g$ does not satisfy triangle inequality in general. 
So we define a function $\rho_g$ from $X\times X$ to $[0,\infty)$ 
by the following:
\begin{align*}
\rho_g(x,y)
=	\inf\{	&
	\rho'_g(x,z_0)+\cdots+\rho'_g(z_i,z_{i+1})+\cdots+\rho'_g(z_{l-1},y)
	\st	\\
&	\quad	l<\omega\text{ and }z_0,\ldots,z_{l-1}\in X
	\}.
\end{align*}

Note that, 
since $f$ is increasing,  
$\rho'_g(x,y)
	\geq f(\max\{c(x),c(y)\})\cdot\abs{c(x)-c(y)}
	\geq\abs{F(c(x))-F(c(y))}$. 
Hence we have 
$\rho_g(x,y)\geq\abs{F(c(x))-F(c(y))}$, because 
\begin{align*}
&	\rho'_g(x,z_0)+\cdots+\rho'_g(z_{l-1},y)	\\
&	\geq\abs{F(c(x))-F(c(z_0))}+\cdots+\abs{F(c(z_{l-1}))-F(c(y))}	\\
&	\geq\abs{F(c(x))-F(c(y))}.
\end{align*}

\begin{claim}\label{claim:magnification}
Let $x,y$ be points of $X$. 
If $x,y\in X\ssm K_{n-1}$, $n<\omega$, 
then 
$\rho_g(x,y)\geq f(\frac{n}{2})\cdot d(x,y)
	\geq g(n)\cdot d(x,y)$. 
\end{claim}

\begin{proof}
We may assume that 
$c(x)=r\geq s=c(y)$, 
$x\in K_m\ssm K_{m-1}$ and 
$y\in K_m$ 
for some $m\geq n$. 
By the definition of $c$, we have $s\geq n$. 
Since $f$ is increasing, 
it suffices to show that 
$\rho'_g(x,z_0)+\cdots+\rho'_g(z_{l-1},y)\geq f(\frac{s}{2})\cdot d(x,y)$ 
for any $l<\omega$, $z_0,\ldots,z_{l-1}\in X$. 

\emph{Case 1.} 
Assume that 
$c(z_i)>\frac{s}{2}$ 
for all $i<l$. 
Since $f$ is increasing, 
the definition of $\rho'_g$ yields 
\begin{align*}
\rho'_g(x,z_0)+\cdots +\rho'_g(z_{l-1},y)
&> f(\tfrac{s}{2})\cdot(\rho(x,z_0)+\cdots+\rho(z_{l-1},y))	\\
&\geq f(\tfrac{s}{2})\cdot\rho(x,y)	\\
&\geq f(\tfrac{s}{2})\cdot d(x,y). 
\end{align*}

\emph{Case 2.}
Assume that $c(z_i)\leq\frac{s}{2}$ for some $i<l$. 
Fix such an $i$ and then 
we have the following: 
\begin{align*}
\rho'_g(x,z_0)+\cdots +\rho'_g(z_{i-1},z_i)
&\geq \rho_g(x,z_i)\geq F(c(x))-F(c(z_i)),	\\
\rho'_g(z_i,z_{i+1})+\cdots +\rho'_g(z_{l-1},y)
&\geq \rho_g(z_i,y)\geq F(c(y))-F(c(z_i)).
\end{align*}
Hence it holds that 
\begin{align*}
\rho'_g(x,z_0)+\cdots +\rho'_g(z_{l-1},y)
&	\geq (F(r)-F(c(z_i)))+(F(s)-F(c(z_i)))	\\
&	\geq (F(r)-F(\tfrac{s}{2}))+(F(s)-F(\tfrac{s}{2}))	\\
&	\geq (r-\tfrac{s}{2})f(\tfrac{s}{2})+\tfrac{s}{2}f(\tfrac{s}{2})	\\
&	=	rf(\tfrac{s}{2}).
\end{align*}
On the other hand, 
$d(x,y)\leq r$, 
because $x,y\in K_m$ 
and $r=c(x)\geq\diam_d K_m$. 
So we have 
\[
\rho'_g(x,z_0)+\cdots +\rho'_g(z_{l-1},y)
\geq f(\tfrac{s}{2})\cdot d(x,y). 
\]

This concludes the proof of the claim.
\end{proof}

Clearly $\rho_g$ is symmetric and satisfies the triangle inequality. 
Since $f(s)\geq 1$ for all $s\in[0,\infty)$, 
Claim~\ref{claim:magnification} implies that $\rho_g$ is a metric on $X$. 
Moreover, $\rho_g$ is proper because $\rho_g\geq\rho$ and $\rho$ is proper. 
It is easy to see that 
$\rho_g$ is compatible with the topology of $(X,d)$. 

Finally, 
we define a metric $d_g$ using $\rho_g$. 
Let $\delta$ be the continuous function from $X$ to $[0,\infty)$ 
which is obtained by applying 
Lemma~\ref{lem:stepfunction} 
to 
$\langle K_n\st n<\omega\rangle$ 
and 
$\langle n^2\st n<\omega\rangle$. 
Note that, 
for $n<\omega$, $x\in K_{n-1}$ and 
$y\in X\ssm K_n$ 
we have 
$\delta(y)-\delta(x)
\geq n$. 
Define $d_g$ by letting 
$d_g(x,y)=\max\{\abs{\delta(x)-\delta(y)},\rho_g(x,y)\}$ 
for $x,y\in X$. 
Then $d_g$ satisfies all requirements of the lemma. 
\end{proof}

\begin{lem}\label{lem:upper}
For any locally compact separable metrizable space $X$, 
we have $\fha(X)\leq\fd$. 
\end{lem}

\begin{proof}
Fix a metric $d$ on $X$, 
and choose 
a sequence $\langle K_n\st n<\omega\rangle$ 
of compact sets of $X$ 
that 
meets the requirement in Lemma~\ref{lem:expand}. 
For each $g\in\omega^\omega$, 
let $d_g$ be the metric on $X$ 
which is obtained by applying Lemma~\ref{lem:expand} 
to $(X,d)$, $\langle K_n\st n<\omega\rangle$ and $g$. 

Choose a subset 
$\mathcal{F}$ of $\omega^\omega$ of size $\fd$ 
which is cofinal with respect to $\domn$. 
We will prove that 
$\stonecech{X}\cptequiv\sup\{\higson{X}{d_g}\st g\in\mathcal{F}\}$. 
It suffices to show that, 
for any two disjoint closed sets $A,B$ of $X$ 
there is a $g\in\mathcal{F}$ such that 
$\cptsep{A}{B}{\higson{X}{d_g}}$. 

For $n<\omega$, 
let $\Delta_n=K_{n+2}\ssm\operatorname{int}K_n$. 
Note that $\Delta_n\subseteq X\ssm K_{n-1}$ for each $n<\omega$. 
Since $A,B$ are disjoint closed sets and each $\Delta_n$ is compact, 
we have $d(A\cap\Delta_n,B\cap\Delta_n)>0$ 
if 
$A\cap\Delta_n\neq\emptyset\neq B\cap\Delta_n\neq\emptyset$. 
Define $h_{A,B}\in\omega^\omega$ 
as follows: 
For $n<\omega$ with 
$A\cap\Delta_n\neq\emptyset\neq B\cap\Delta_n\neq\emptyset$, 
let 
\[	h_{A,B}(n)
	=\left\lceil\frac{n}{d(A\cap\Delta_n,B\cap\Delta_n)}\right\rceil,	\] 
(where $\lceil r\rceil$ denotes the smallest integer not smaller than $r$)
and 
otherwise $h_{A,B}(n)$ is arbitrary. 
Find $g\in\mathcal{F}$ and $N<\omega$ 
such that $h_{A,B}(n)\leq g(n)$ for $n>N$. 

We claim that, 
for every $M\geq N$ and $x\in X\ssm K_{M+1}$ 
we have $d_g(x,A)+d_g(x,B)\geq M$, 
and hence $\cptsep{A}{B}{\higson{X}{d_g}}$. 
Fix $M<\omega$ and $x\in X\ssm K_{M+1}$. 
Since $d_g$ is a proper metric, 
we can find $a\in A$ and $b\in B$ such that 
$d_g(x,A)+d_g(x,B)=d_g(x,a)+d_g(x,b)$ holds. 
Choose $n_a,n_b<\omega$ so that 
$a\in K_{n_a}\ssm K_{n_a-1}$ 
and 
$b\in K_{n_b}\ssm K_{n_b-1}$, 
and let $n=\min\{n_a,n_b\}$. 

\emph{Case 1.}
$n\leq M$. 
Since $x\in X\ssm K_{M+1}$ and $d_g(K_M,X\ssm K_{M+1})\geq M$, 
we have $d_g(a,x)\geq M$ or $d_g(b,x)\geq M$. 

\emph{Case 2.}
$n>M$. 
By the triangle inequality, 
it suffices to show that $d_g(a,b)\geq M$. 
If $\abs{n_a-n_b}\leq 1$, 
then $a,b\in \Delta_{n-1}$, 
and hence 
we have 
\begin{align*}
d_g(a,b)
	&\geq g(n-1)\cdot d(a,b)	\\	
	&\geq h_{A,B}(n-1)\cdot d(a,b)	\\
	&\geq h_{A,B}(n-1)\cdot d(A\cap\Delta_{n-1},B\cap\Delta_{n-1})	\\
	&\geq n-1 \geq M. 
\end{align*}
Otherwise, 
we have 
$d_g(a,b)\geq d_g(K_n,X\ssm K_{n+1})\geq n> M$.

This concludes the proof. 
\end{proof}

Now we have the following theorem. 

\begin{thm}
Let $X$ be a locally compact separable metrizable space. 
If 
$\deriv{X}$ 
is not compact, 
then 
$\fsa(X)=\fha(X)=\fd$, 
and otherwise $\fsa(X)=\fha(X)=1$. 
\end{thm}

\section{It may be further than $\fd$ miles}\label{sec:morethand}

The cardinal $\fha(X)$ is defined 
for locally compact separable metrizable spaces $X$, 
while $\fsa(X)$ is defined for any metrizable space $X$. 
By Theorem~\ref{thm:saone} and Lemma~\ref{lem:lower}, 
either $\fsa(X)\geq\fd$ or $\fsa(X)=1$ holds for any $X$. 
In this section, 
we show the existence of a metrizable space $X$ 
for which $\fsa(X)>\fd$ holds. 

For a topological space $X$, 
$e(X)$, the \emph{extent} of $X$, 
is defined by 
$e(X)=\sup\{\size{D}\st 
D\subseteq X\text{ and }D\text{ is closed discrete}\}+\aleph_0$. 

\begin{defn}
For an infinite cardinal $\kappa$, 
define $\log\kappa$ by letting 
$\log\kappa=\min\{\theta\st 2^\theta\geq\kappa\}$. 
\end{defn}

It is easy to see that, 
for a set $C$ of infinite cardinals, 
we have $\log(\sup C)=\sup\{\log\kappa\st\kappa\in C\}$.

\begin{prop}
Let $X$ be a metrizable space. 
If 
$\deriv{X}$ 
is not compact, 
then 
$\fsa(X)\geq\log e(\deriv{X})$. 
\end{prop}

\begin{proof}
It suffices to show that, 
for infinite cardinals $\kappa$ and $\lambda$, 
if 
$\deriv{X}$ 
has a closed discrete subset of size $\kappa$ 
and $\lambda=\log\kappa$, 
then $\fsa(X)\geq\lambda$. 

Suppose that $D$ is a set of compatible metrics on $X$ 
and $\size{D}=\mu<\lambda$. 
We will show that 
$\stonecech{X}\not\cptequiv\sup\{\smirnov{X}{\rho}\st\rho\in D\}$. 
Since we have $\fsa(X)\geq\fd$ by Lemma~\ref{lem:lower}, 
we 
may 
assume that $\mu\geq\fd$. 

Choose a subset $H$ of $\omega^\omega$ of size $\fd$ 
which is cofinal with respect to $\leq$. 

Fix a closed discrete subset 
$A=\{a_\xi\st\xi<\kappa\}$ of 
$\deriv{X}$. 
As in the proof of Lemma~\ref{lem:lower}, 
we choose a neighborhood $U_\xi$ of $a_\xi$ 
and a sequence $\langle b_{\xi,i}\st i<\omega\rangle$ 
in $U_\xi\ssm\{a_\xi\}$
for $\xi<\kappa$ 
so that, 
\begin{enumerate}
\item for each $\xi<\kappa$, 
	$\langle b_{\xi,i}\st i<\omega\rangle$ converges to $a_\xi$, 
\item if $\xi<\eta<\kappa$ then $U_\xi\cap U_\eta=\emptyset$, and 
\item for any $\varphi\in\omega^\kappa$, 
	the set $\{b_{\xi,\varphi(\xi)}\st\xi<\kappa\}$ 
	has no accumulating point. 
\end{enumerate}

For each $\rho\in D$ and $\xi<\kappa$, 
define $g^\xi_\rho\in\omega^\omega$ 
by letting 
\[
g^\xi_\rho(m)=
	\min\left\{k<\omega\st 
	\forall i\geq k\,
	\left(\rho(a_\xi,b_{\xi,i})<\frac{1}{m+1}\right)\right\}
\]
for $m<\omega$, 
and choose $h^\xi_\rho\in H$ 
so that $g^\xi_\rho\leq h^\xi_\rho$. 

Since 
$\fd\leq\mu=\size{D}
<\lambda=\log\kappa$, 
we have $\fd^\mu=2^\mu<\kappa$, 
and hence there are $K\in[\kappa]^\kappa$ 
and $\{h^\xi\st\xi\in K\}$ 
such that, 
for each $\xi\in K$, 
$h^\xi_\rho=h^\xi$ for all $\rho\in D$. 

Fix a countable set $\{\xi_n\st n<\omega\}\subseteq K$. 
Let $b_n=b_{\xi_n,h^{\xi_n}(n)}$ 
and $B=\{b_n\st n<\omega\}$. 
By the choice of $A$, $U_\xi$'s and $b_{\xi,i}$'s, 
$A\cap B=\emptyset$ and $B$ is closed in $X$. 
Also, by the choice of $h^\xi_\rho$'s, 
for each $\rho\in D$ and $n<\omega$ 
we have $\rho(a_{\xi_n},b_n)\leq\frac{1}{n+1}$. 

Now it is easy to see that 
$\cptnonsep{A}{B}{\sup\{\smirnov{X}{\rho}\st\rho\in D\}}$, 
and hence 
$\stonecech{X}\not\cptequiv\sup\{\smirnov{X}{\rho}\st\rho\in D\}$. 
\end{proof}

\begin{cor}
Let $X_\kappa=\kappa\times(\omega+1)$, 
where $\kappa$ is equipped with the discrete topology 
and $\omega+1$ is equipped with the usual order topology. 
If $\kappa>2^\fd$, then $\fsa(X_\kappa)>\fd$. 
\end{cor}

%


\begin{thebibliography}{1}

\bibitem{BaY:ring}
B.~J. Ball and S.~Yokura.
\newblock Compactifications determined by subsets of {$C^*(X)$}.
\newblock {\em Topology Appl.}, 13:1--13, 1982.

\bibitem{Chan:cpt}
R.~E. Chandler.
\newblock {\em Hausdorff Compactifications}.
\newblock Marcel Dekker Inc., New York, 1976.

\bibitem{DKU:higsoncorona}
A.~N. Dranishnikov, J.~Keesling, and V.~V. Uspenskij.
\newblock On the {Higson} corona of uniformly contractible spaces.
\newblock {\em Topology}, 37:791--803, 1998.

\bibitem{GJ:rings}
L.~Gillman and M.~Jerison.
\newblock {\em Rings of continuous functions}.
\newblock Van Nostrand, 1960.

\bibitem{KTY:babylon}
M.~Kada, K.~Tomoyasu, and Y.~Yoshinobu.
\newblock How many miles to $\beta\omega$? --- {Approximating} $\beta\omega$ by
  metric-dependent compactifications.
\newblock {\em Topology Appl.}, 145:277--292, 2004.

\bibitem{KT:approx}
K.~Kawamura and K.~Tomoyasu.
\newblock Approximations of {Stone--\v{C}ech} compactifications by {Higson}
  compactifications.
\newblock {\em Colloquium Mathematicum}, 88:75--92, 2001.

\bibitem{Wo:unifcpt}
R.~G. Woods.
\newblock The minimum uniform compactification of a metric space.
\newblock {\em Fund. Math.}, 147:39--59, 1995.

\end{thebibliography}

\vspace*{12pt}
\indent
Masaru Kada,
    Information Processing Center, 
	Kitami Institute of Technology. 
	Kitami 090--8507 JAPAN.\par
%
%
\indent 
\textit{Current address}:
	Graduate School of Science, 
	Osaka Prefecture University. 
	Sakai, Osaka 599--8531 JAPAN.\par
\indent
%
%
e-mail: kada@mi.s.osakafu-u.ac.jp\par
\indent
Kazuo Tomoyasu,
	General Education, Miyakonojo National College of Technology,
	Miyakonojo-shi, Miyazaki 885--8567 JAPAN.\par
\indent
e-mail:tomoyasu@cc.miyakonojo-nct.ac.jp\par
\indent
Yasuo Yoshinobu
	Graduate School of Information Science,	
	Nagoya University. 
	Nagoya 464--8601 JAPAN.\par
\indent
e-mail:yosinobu@math.nagoya-u.ac.jp

\end{document}